\documentclass[final,leqno,onefignum,onetabnum]{siamltex1213}

\usepackage{graphicx,epsfig,color}
\usepackage{booktabs}
\usepackage[notref,notcite]{showkeys}
\usepackage{amsmath,cases}
\usepackage{amssymb,bm}
\usepackage{multirow}
\newtheorem{remark}[theorem]{Remark}
\input{siam.sty}

\title{Refined upper bounds for the convergence of the randomized extended Kaczmarz and Gauss-Seidel algorithms}

\author{Kui Du\thanks{School of Mathematical Sciences and Fujian Provincial Key Laboratory of Mathematical Modeling and High-Performance Scientific Computing, Xiamen University, Xiamen 361005, China ({\tt kuidu@xmu.edu.cn}). The work of the author was supported by the National Natural Science Foundation of China (No.11771364 and No.91430213) and the Fundamental Research Funds for the Central Universities (No.20720160002).}}

\begin{document}
\maketitle

\pagestyle{myheadings} \thispagestyle{plain} \markboth{KUI DU}{RANDOMIZED EXTENDED KACZMARZ AND GAUSS-SEIDEL}

\begin{abstract} The randomized extended Kaczmarz and Gauss-Seidel algorithms have attracted much attention because of their ability to treat all types of linear systems (consistent or inconsistent, full rank or rank-deficient). 
In this paper, we interpret the randomized extended Kaczmarz and Gauss-Seidel algorithms  as specific combinations of the randomized Kaczmarz and Gauss-Seidel algorithms and present refined upper bounds for their convergence. 
\end{abstract}

\begin{keywords} randomized extended Kaczmarz algorithm, randomized extended Gauss-Seidel algorithm, Moore-Penrose pseudoinverse solution, convergence analysis
\end{keywords}

\begin{AMS} 65F10, 65F20
\end{AMS} 

\section{Introduction} Due to the better performance in many situations than existing classical iterative algorithms, randomized iterative algorithms for solving a linear system of equations 
\beq\label{lin} {\bf Ax=b},\quad \mbf A\in\mbbr^{m\times n},\quad \mbf b\in\mbbr^m,\eeq have attracted much attention recently; see, for example, \cite{strohmer2009rando,leventhal2010rando,needell2010rando,eldar2011accel,zouzias2013rando,needell2014paved,dumitrescu2015relat,ma2015conve,gower2015rando,liu2016accel,needell2016stoch,hefny2017rows} and the references therein. In this paper, we consider the randomized Kaczmarz (RK) algorithm \cite{strohmer2009rando}, the randomized Gauss-Seidel (RGS) algorithm \cite{leventhal2010rando}, the randomized extended Kaczmarz (REK) algorithm \cite{zouzias2013rando}, and the randomized extended Gauss-Seidel (REGS) algorithm \cite{ma2015conve}. 
Let $\mbf A^\dag$ denote the Moore-Penrose pseudoinverse \cite{ben2003gener} of $\bf A$. We summarize the convergence of RK, RGS, REK, and REGS  in expectation to the Moore-Penrose pseudoinverse solution $\mbf A^\dag\mbf b$ for all types of linear systems in Table \ref{t1}. 

\begin{table}[htp]
\caption{Summary of the convergence of {\rm RK}, {\rm RGS}, {\rm REK}, and {\rm REGS} in expectation to the Moore-Penrose pseudoinverse solution $\mbf A^\dag\mbf b$ for all types of linear systems: {\rm Y} means the algorithm is convergent and {\rm N} means not.}  \label{t1}
\begin{center} 
\begin{tabular}{c|c|c|c|c|c} \toprule
linear system (\ref{lin}) & $\rank(\mbf A)$ & RK & RGS &REK & REGS \\ 
\hline
consistent & $=n$ & Y   &  Y  & Y  & Y  \\ 
 consistent & $<n$ & Y &  N   & Y  & Y    \\
 inconsistent & $=n$  &  N & Y   & Y  & Y  \\
 inconsistent & $<n$ & N  &  N   & Y   & Y \\ \bottomrule
\end{tabular}
\end{center}
\end{table}

{\it Main contributions}. We show that REK is essentially an RK-RK approach (see Remark \ref{RKRK}) and that REGS is essentially an RGS-RK approach (see Remark \ref{RGSRK}). We present refined upper bounds for the convergence of REK and REGS. These bounds hold for all types of linear systems (consistent or inconsistent, overdetermined or underdetermined, $\mbf A$ has full column rank or not) and are attainable. In addition, we point out that the proof for Theorem 4.1 of \cite{ma2015conve} is incomplete and we resolve this issue. 


{\it Organization of the paper}. In the rest of this section, we give some notation and preliminaries. In section 2, we review the randomized Kaczmarz algorithm and the randomized extended Kaczmarz algorithm. We present a slightly different variant of REK and prove its convergence. In section 3, we review the randomized Gauss-Seidel algorithm and the randomized extended Gauss-Seidel algorithm. We show that the convergence analysis for REGS of \cite{ma2015conve} is incomplete. We present a mathematically equivalent variant of REGS and prove its convergence. Numerical examples are given in section 4 to illustrate the theoretical results. We present brief concluding remarks in section 5.

{\it Notation and preliminaries}. For any random variable $\bm\xi$, let $\mbbe\bem\bm\xi\eem$ denote its expectation. For an integer $m\geq 1$, let $[m]:=\{1,2,3,\ldots,m\}$. Throughout the paper all vectors are assumed to be column vectors. For any vector $\mbf u\in\mbbr^m$, we use $\bf u^\rmt$,  $u_i$, and $\|\mbf u\|_2$ to denote the transpose, the $i$th entry, and the Euclidean norm of $\mbf u$, respectively. We use ${\bf e}_j$ to denote the $j$th column of the identity matrix $\mbf I$ whose order is clear from the context. For any matrix $\mbf A\in\mbbr^{m\times n}$, we use $\mbf A^\rmt$, $\|\mbf A\|_\rmf$, $\rank(\mbf A)$, $\ran(\mbf A)$, $\nul(\mbf A)$, $\sigma_{1}(\mbf A)$, and $\sigma_{r}(\mbf A)$ to denote the transpose, the Frobenius norm, the rank, the column space, the nullspace, the largest singular value, and the smallest nonzero singular value of $\mbf A$, respectively. 
We denote the columns and rows of $\bf A$  by $\{\mbf a_j\}_{j=1}^n$ and $\{\wt{\mbf a}_i^\rmt\}_{i=1}^m$, respectively. That is to say, $${\bf A}=\bem {\bf a}_1 & {\bf a}_2 &\cdots& {\bf a}_n\eem,\quad {\bf A}^\rmt=\bem\wt{\bf a}_1 & \wt{\bf a}_2 & \cdots & \wt{\bf a}_m \eem.$$  All the convergence results depend on the positive number $\rho$ defined as $$\rho:=1-\frac{\sigma_{r}^2(\bf A)}{\|\mbf A\|_\rmf^2}.$$ The following lemmas will be used extensively in this paper. Their proofs are straightforward. 

\lemma\label{leq} Let $\mbf A$ be any nonzero real matrix. For every $\mbf u\in\ran(\mbf A)$, it holds \beqs\label{est1}\mbf u^\rmt\l(\mbf I-\frac{\bf AA^\rmt}{\|\mbf A\|_\rmf^2}\r)\mbf u\leq\rho\|\mbf u\|_2^2.\eeqs The equality holds if $\sigma_1(\mbf A)=\sigma_r(\mbf A)$, i.e., all the nonzero singular values of $\mbf A$ are the same. 
\endlemma
\lemma\label{proj} Let $\mbf a$ be any nonzero vector. Then $$\l(\frac{\mbf a\mbf a^\rmt}{\|\mbf a\|_2^2}\r)^2=\frac{\mbf a\mbf a^\rmt}{\|\mbf a\|_2^2},\quad \l(\mbf I-\frac{\mbf a\mbf a^\rmt}{\|\mbf a\|_2^2}\r)^2=\mbf I-\frac{\mbf a\mbf a^\rmt}{\|\mbf a\|_2^2}.$$
\endlemma

\section{Randomized Kaczmarz and its extension} 

Strohmer and Vershynin \cite{strohmer2009rando} proposed the following randomized Kaczmarz algorithm (Algorithm 1). 
\begin{table}[htp]
\begin{center}
\begin{tabular*}{115mm}{l}
\toprule {\bf Algorithm 1}. Randomized Kaczmarz \cite{strohmer2009rando} for $\bf Ax=b$
\\ \hline \vspace{-3mm}\\ 
\qquad Initialize $\mbf x^0\in\mbbr^n$ \\
\qquad {\bf for} $k=1,2,\ldots$ {\bf do}\\
\qquad\qquad Pick $i\in[m]$ with probability $\|\wt{\bf a}_i\|_2^2/\|{\bf A}\|_\rmf^2$\\
\qquad\qquad Set $\dsp{\bf x}^k={\bf x}^{k-1}-\frac{\wt{\bf a}_i^\rmt{\bf x}^{k-1}-b_i}{\|\wt{\bf a}_i\|_2^2}\wt{\bf a}_i$\\
\bottomrule
\end{tabular*}
\end{center}
\end{table}

If $\bf Ax=b$ is consistent, Zouzias and Freris \cite[Theorem 3.4]{zouzias2013rando} proved that RK with initial guess $\mbf x^0\in\ran(\mbf A^\rmt)$ generates $\mbf x^k$ which  converges linearly in expectation to the Moore-Penrose pseudoinverse solution ${\bf A^\dag b}$: $$\mbbe\bem\|{\bf x}^k-{\bf A^\dag b} \|_2^2\eem\leq\rho^k\|{\bf x}^0-{\bf A^\dag b}\|_2^2.$$ By the same approach as used in the proof of Theorem 3.2 of \cite{zouzias2013rando}, we can prove the following theorem, which will be used to prove the refined upper bound for the convergence of REK.

\theorem\label{estzk} Let ${\bf A}\in\mbbr^{m\times n}$ and ${\bf b}\in\mbbr^m$. Let ${\bf z}^k$ denote the $k$th iterate of {\rm RK} applied to ${\bf A^\rmt z=0}$ with initial guess ${\bf z}^0\in{\bf b}+\ran({\bf A})$. In exact arithmetic, it holds \beqs\label{bound1}\mbbe\bem\|{\bf z}^k-{(\bf I-AA^\dag)b}\|_2^2\eem\leq\rho^k\|\mbf z^0- {(\mbf I-\bf AA^\dag) b} \|_2^2.\eeqs
\endtheorem

\proof
The iteration is $${\bf z}^k={\bf z}^{k-1}-\frac{{\bf a}_j^\rmt{\bf z}^{k-1}}{\|{\bf a}_j\|_2^2}{\bf a}_j.$$  By $\mbf a_j^\rmt{\bf(I- AA^\dag) b}=0$ (since ${\bf A^\rmt (I- AA^\dag) b=0}$), we have \beqas
{\bf z}^k-{\bf (I- AA^\dag) b}&=&{\bf z}^{k-1}-{\bf (I- AA^\dag) b}-\frac{{\bf a}_j^\rmt{\bf z}^{k-1}-\mbf a_j^\rmt{\bf(I- AA^\dag) b}}{\|{\bf a}_j\|_2^2}{\bf a}_j\\ \nn &=&{\bf z}^{k-1}-{\bf (I- AA^\dag) b}-\frac{{\bf a}_j^\rmt({\bf z}^{k-1}-{\bf(I- AA^\dag) b})}{\|{\bf a}_j\|_2^2}{\bf a}_j\\\nn&=&\l({\bf I}-\frac{{\bf a}_j{\bf a}_j^\rmt}{\|{\bf a}_j\|_2^2}\r)({\bf z}^{k-1}-{\bf(I- AA^\dag) b}).\eeqas 
By ${\bf z}^0\in{\bf b}+\ran({\bf A})$ and ${\bf AA^\dag b}\in\ran({\bf A})$, we have ${\bf z}^0-{\bf(I- AA^\dag) b}\in\ran(\bf A)$. Then it is easy to show that ${\bf z}^k-{\bf(I- AA^\dag) b}\in\ran(\bf A)$ by induction. Let $\mbbe_{k-1}\bem\cdot\eem$ denote the conditional expectation conditioned on the first $k-1$ iterations of RK. It follows that \beqas&&\mbbe_{k-1}\bem\|{\bf z}^k-{\bf(I- AA^\dag) b}\|_2^2\eem\\&=&\mbbe_{k-1}\bem ({\bf z}^{k}-{\bf(I- AA^\dag) b})^\rmt({\bf z}^{k}-{\bf(I- AA^\dag) b})\eem\\&=&\mbbe_{k-1}\bem({\bf z}^{k-1}-{\bf(I- AA^\dag) b})^\rmt\dsp\l({\bf I}-\frac{{\bf a}_j{\bf a}_j^\rmt}{\|{\bf a}_j\|_2^2}\r)^2({\bf z}^{k-1}-{\bf(I- AA^\dag) b})\eem
\\&=&\mbbe_{k-1}\bem({\bf z}^{k-1}-{\bf(I- AA^\dag) b})^\rmt\dsp\l({\bf I}-\frac{{\bf a}_j{\bf a}_j^\rmt}{\|{\bf a}_j\|_2^2}\r)({\bf z}^{k-1}-{\bf(I- AA^\dag) b})\eem 
\\ &=&({\bf z}^{k-1}-{\bf(I- AA^\dag) b})^\rmt\l({\bf I}-\frac{\bf AA^\rmt}{\|{\bf A}\|_\rmf^2}\r)({\bf z}^{k-1}-{\bf(I- AA^\dag) b})\\&\leq&\rho\|{\bf z}^{k-1}-{\bf(I- AA^\dag) b}\|_2^2. \quad (\mbox{by Lemma \ref{leq}})\eeqas Taking expectation gives $$\mbbe\bem\|{\bf z}^k-{\bf(I- AA^\dag) b}\|_2^2\eem\leq \rho\mbbe\bem\|{\bf z}^{k-1}-{\bf(I- AA^\dag) b}\|_2^2\eem.$$ Unrolling the recurrence yields the result. \endproof

If $\bf Ax=b$ is inconsistent, Needell \cite{needell2010rando} and Zouzias and Freris \cite{zouzias2013rando} showed that RK does not converge to $\mbf A^\dag\mbf b$. To resolve this problem, Zouzias and Freris \cite{zouzias2013rando} proposed the following randomized extended Kaczmarz algorithm (here we call it REK-ZF, see Algorithm 2). They proved the convergence bound \beq\label{rekb} \mbbe\bem\|{\bf x}^k-\mbf A^\dag\mbf b\|_2^2\eem\leq\rho^{\lf k/2\rf}(1+2\sigma_{1}^2(\mbf A)/\sigma_{r}^2(\mbf A))\|\mbf A^\dag\mbf b\|_2^2.\eeq  
\begin{table}[htp]
\begin{center}
\begin{tabular*}{115mm}{l}
\toprule {\bf Algorithm 2}. REK-ZF \cite{zouzias2013rando}
\\ \hline \vspace{-3mm}\\ 
\qquad Initialize $\mbf x^0\in\ran({\bf A}^\rmt)$ and ${\bf z}^0\in\mbf b+\ran(\mbf A)$ \\
\qquad {\bf for} $k=1,2,\ldots$ {\bf do}\\
\qquad\qquad Pick $j\in[n]$ with probability $\|{\bf a}_j\|_2^2/\|{\bf A}\|_\rmf^2$\\
\qquad\qquad Set $\dsp{\bf z}^k={\bf z}^{k-1}-\frac{{\bf a}_j^\rmt{\bf z}^{k-1}}{\|{\bf a}_j\|_2^2}{\bf a}_j$\\
\qquad\qquad Pick $i\in[m]$ with probability $\|\wt{\bf a}_i\|_2^2/\|{\bf A}\|_\rmf^2$\\
\qquad\qquad Set $\dsp{\bf x}^k={\bf x}^{k-1}-\frac{\wt{\bf a}_i^\rmt{\bf x}^{k-1}-b_i+z_i^{k-1}}{\|\wt{\bf a}_i\|_2^2}\wt{\bf a}_i$\\
\bottomrule
\end{tabular*}
\end{center}
\end{table}

\remark\label{RKRK} The original randomized extended Kaczmarz algorithm of {\rm\cite{zouzias2013rando}} uses $\mbf x^0=\mbf 0$ and ${\bf z}^0=\mbf b$. Essentially, {\rm REK-ZF} is an {\rm RK-RK} approach: $\mbf z^k$ is the $k$th iterate of {\rm RK} applied to $\bf A^\rmt z=0$ with initial guess $\mbf z^0$, and $\mbf x^k$ is one step {\rm RK} update for the linear system ${\bf Ax=b-z}^{k-1}$ from $\mbf x^{k-1}$.
\endremark

Next we present a slightly different randomized extended Kaczmarz algorithm (we call it REK-S, see Algorithm 3) which generates $\mbf x^k$ by one step RK update for the linear system ${\bf Ax=b-z}^k$ (used in REK-S) instead of ${\bf Ax=b-z}^{k-1}$ (used in REK-ZF) from $\mbf x^{k-1}$. 
\begin{table}[htp]
\begin{center}
\begin{tabular*}{115mm}{l}
\toprule {\bf Algorithm 3}. REK-S
\\ \hline \vspace{-3mm}\\ 
\qquad Initialize $\mbf x^0\in\ran(\mbf A^\rmt)$ and ${\bf z}^0\in{\bf b}+\ran({\bf A})$ \\
\qquad {\bf for} $k=1,2,\ldots$ {\bf do}\\
\qquad\qquad Pick $j\in[n]$ with probability $\|{\bf a}_j\|_2^2/\|{\bf A}\|_\rmf^2$\\
\qquad\qquad Set $\dsp{\bf z}^k={\bf z}^{k-1}-\frac{{\bf a}_j^\rmt{\bf z}^{k-1}}{\|{\bf a}_j\|_2^2}{\bf a}_j$\\
\qquad\qquad Pick $i\in[m]$ with probability $\|\wt{\bf a}_i\|_2^2/\|{\bf A}\|_\rmf^2$\\
\qquad\qquad Set $\dsp{\bf x}^k={\bf x}^{k-1}-\frac{\wt{\bf a}_i^\rmt{\bf x}^{k-1}-b_i+z_i^k}{\|\wt{\bf a}_i\|_2^2}\wt{\bf a}_i$\\
\bottomrule
\end{tabular*}
\end{center}
\end{table}
In the following theorem, we present the convergence bound for REK-S, which is obviously better than the bound (\ref{rekb}). Actually, our bound is attainable (see Remark \ref{sharp1}). We note that a refined convergence bound for REK-ZF can be obtained by the same approach.

\theorem\label{main1} Let ${\bf A}\in\mbbr^{m\times n}$ and ${\bf b}\in\mbbr^m$. Let ${\bf x}^k$ denote the $k$th iterate of {\rm REK-S} with $\mbf x^0\in\ran(\mbf A^\rmt)$ and  $\mbf z^0\in\mbf b+\ran(\mbf A)$. In exact arithmetic, it holds \beq\label{bound2}\mbbe\bem\|{\bf x}^k-{\bf A^\dag b}\|_2^2\eem\leq\rho^k\|{\bf x}^0-{\bf A^\dag b}\|_2^2+\frac{\rho^k(1-\rho^k)}{\sigma_{r}^2(\mbf A)}\|{\bf z}^0-{\bf (I-AA^\dag) b}\|_2^2.\eeq
\endtheorem
\proof Let $$\wh{\mbf x}^k={\bf x}^{k-1}-\frac{\wt{\bf a}_i^\rmt{\bf x}^{k-1}-b_i+\mbf e_i^\rmt{\bf (I-AA^\dag) b}}{\|\wt{\bf a}_i\|_2^2}\wt{\bf a}_i.$$ We have \beqas\wh{\mbf x}^k-{\bf A^\dag b}&=&{\bf x}^{k-1}-{\bf A^\dag b}-\frac{\wt{\bf a}_i^\rmt{\bf x}^{k-1}-\mbf e_i^\rmt{\bf AA^\dag b}}{\|\wt{\bf a}_i\|_2^2}\wt{\bf a}_i\\&=&{\bf x}^{k-1}-{\bf A^\dag b}-\frac{\wt{\bf a}_i^\rmt{\bf x}^{k-1}-\wt{\bf a}_i^\rmt{\bf A^\dag b}}{\|\wt{\bf a}_i\|_2^2}\wt{\bf a}_i\\&=&\l(\mbf I-\frac{\wt{\bf a}_i\wt{\bf a}_i^\rmt}{\|\wt{\bf a}_i\|_2^2}\r)({\bf x}^{k-1}-{\bf A^\dag b})\eeqas and $${\mbf x}^k-\wh{\mbf x}^k=\frac{\mbf e_i^\rmt({\bf (I-AA^\dag) b}-\mbf z^k)}{\|\wt{\bf a}_i\|_2^2}\wt{\bf a}_i.$$ By the orthogonality $(\wh{\mbf x}^k-{\bf A^\dag b})^\rmt({\mbf x}^k-\wh{\mbf x}^k)=0$ (which is obvious from the above two equations), we have \beq\label{ksum}\|{\mbf x}^k-{\bf A^\dag b}\|_2^2=\|{\mbf x}^k-\wh{\mbf x}^k\|_2^2+\|\wh{\mbf x}^k-{\bf A^\dag b}\|_2^2.\eeq Let $\mbbe_{k-1}\bem\cdot\eem$ denote the conditional expectation conditioned on the first $k-1$ iterations of REK-S. That is, $$\mbbe_{k-1}\bem\cdot\eem=\mbbe\bem\cdot|j_1,i_1,j_2,i_2,\ldots,j_{k-1},i_{k-1}\eem,$$ where $j_l$ is the $l$th column chosen and $i_l$ is the $l$th row chosen. We denote the conditional expectation conditioned on the first $k-1$ iterations and the $k$th column chosen as $$\mbbe_{k-1}^i\bem\cdot\eem=\mbbe\bem\cdot|j_1,i_1,j_2,i_2,\ldots,j_{k-1},i_{k-1},j_k\eem.$$ Similarly, we denote the conditional expectation conditioned on the first $k-1$ iterations and the $k$th row chosen as $$\mbbe_{k-1}^j\bem\cdot\eem=\mbbe\bem\cdot|j_1,i_1,j_2,i_2,\ldots,j_{k-1},i_{k-1},i_k\eem.$$ Then by the law of total expectation we have $$\mbbe_{k-1}\bem\cdot\eem=\mbbe_{k-1}^j\bem\mbbe_{k-1}^i\bem\cdot\eem\eem.$$
It follows from \beqas\mbbe_{k-1}\bem\|{\mbf x}^k-\wh{\mbf x}^k\|_2^2\eem&=&\mbbe_{k-1}\bem\dsp\frac{(\mbf e_i^\rmt({\bf (I-AA^\dag) b}-\mbf z^k))^2}{\|\wt{\bf a}_i\|_2^2}\eem\\ &=&\mbbe_{k-1}^j\bem\mbbe_{k-1}^i\bem\dsp\frac{(\mbf e_i^\rmt({\bf (I-AA^\dag) b}-\mbf z^k))^2}{\|\wt{\bf a}_i\|_2^2}\eem\eem\\&=&\mbbe_{k-1}^j\bem\dsp\frac{\|\mbf z^k-{\bf (I-AA^\dag) b}\|_2^2}{\|\mbf A\|_\rmf^2}\eem\\&=&\frac{1}{\|\mbf A\|_\rmf^2}\mbbe_{k-1}\bem\|\mbf z^k-{\bf (I-AA^\dag) b}\|_2^2\eem\eeqas that \beqa\nn\mbbe\bem\|{\mbf x}^k-\wh{\mbf x}^k\|_2^2\eem&=&\frac{1}{\|\mbf A\|_\rmf^2}\mbbe\bem\|\mbf z^k-{\bf (I-AA^\dag) b}\|_2^2\eem\\\label{ksum1}&\leq&\frac{\rho^k}{\|\mbf A\|_\rmf^2}\|\mbf z^0-{\bf (I-AA^\dag) b}\|_2^2. \quad (\mbox{by Theorem \ref{estzk}})\eeqa By $\mbf x^0 \in\ran(\mbf A^\rmt)$ and ${\bf A^\dag b}\in\ran(\mbf A^\rmt)$, we have $\mbf x^0-{\bf A^\dag b}\in\ran(\mbf A^\rmt)$. Then it is easy to show that $\mbf x^k-{\bf A^\dag b}\in\ran(\mbf A^\rmt)$ by induction. It follows from \beqas\mbbe_{k-1}\bem\|\wh{\mbf x}^k-{\bf A^\dag b}\|_2^2\eem&=&\mbbe_{k-1}\bem(\wh{\mbf x}^k-{\bf A^\dag b})^\rmt(\wh{\mbf x}^k-{\bf A^\dag b})\eem\\&=&\mbbe_{k-1}\bem\dsp({\bf x}^{k-1}-{\bf A^\dag b})^\rmt\l(\mbf I-\frac{\wt{\bf a}_i\wt{\bf a}_i^\rmt}{\|\wt{\bf a}_i\|_2^2}\r)^2({\bf x}^{k-1}-{\bf A^\dag b})\eem\\&=&\mbbe_{k-1}\bem\dsp({\bf x}^{k-1}-{\bf A^\dag b})^\rmt\l(\mbf I-\frac{\wt{\bf a}_i\wt{\bf a}_i^\rmt}{\|\wt{\bf a}_i\|_2^2}\r)({\bf x}^{k-1}-{\bf A^\dag b})\eem\\&=&\dsp({\bf x}^{k-1}-{\bf A^\dag b})^\rmt\l(\mbf I-\frac{\bf A^\rmt A}{\|{\bf A}\|_\rmf^2}\r)({\bf x}^{k-1}-{\bf A^\dag b})\\&\leq&\rho\|{\bf x}^{k-1}-{\bf A^\dag b}\|_2^2\quad (\mbox{by Lemma \ref{leq}})\eeqas that \beq\label{ksum2}\mbbe\bem\|\wh{\mbf x}^k-{\bf A^\dag b}\|_2^2\eem\leq\rho\mbbe\bem\|{\bf x}^{k-1}-{\bf A^\dag b}\|_2^2\eem.\eeq Combining (\ref{ksum}), (\ref{ksum1}), and (\ref{ksum2}) yields \beqas\mbbe\bem[\|{\mbf x}^k-{\bf A^\dag b}\|_2^2\eem&=&\mbbe\bem\|{\mbf x}^k-\wh{\mbf x}^k\|_2^2\eem+\mbbe\bem\|\wh{\mbf x}^k-{\bf A^\dag b}\|_2^2\eem\\&\leq&\frac{\rho^k}{\|\mbf A\|_\rmf^2}\|\mbf z^0-{\bf (I-AA^\dag) b}\|_2^2+\rho\mbbe\bem\|{\bf x}^{k-1}-{\bf A^\dag b}\|_2^2\eem\\&\leq&\cdots\leq\rho^k\|{\bf x}^0-{\bf A^\dag b}\|_2^2+\frac{\rho^k}{\|\mbf A\|_\rmf^2}\|\mbf z^0-{\bf (I-AA^\dag) b}\|_2^2\sum_{l=0}^{k-1}\rho^l\\&=&\rho^k\|{\bf x}^0-{\bf A^\dag b}\|_2^2+\frac{\rho^k}{\|\mbf A\|_\rmf^2}\|\mbf z^0-{\bf (I-AA^\dag) b}\|_2^2\frac{1-\rho^k}{1-\rho}\\&=&\rho^k\|{\bf x}^0-{\bf A^\dag b}\|_2^2+\frac{\rho^k(1-\rho^k)}{\sigma_{r}^2(\mbf A)}\|\mbf z^0-{\bf (I-AA^\dag) b}\|_2^2.\eeqas 
This completes the proof. \endproof

\remark The vector $\wh{\mbf x}^k$ used in the proof is actually one step {\rm RK} update for the linear system ${\bf Ax=AA^\dag b}$ from $\mbf x^{k-1}$.
\endremark
\remark\label{sharp1} By Lemma {\rm 1}, if $\sigma_1(\mbf A)=\sigma_r(\mbf A)$, then all the inequalities in the proofs of Theorems {\rm \ref{estzk}} and {\rm \ref{main1}} become equalities. 
\endremark

\section{Randomized Gauss-Seidel and its extension}
Leventhal and Lewis \cite{leventhal2010rando} proposed the following randomized Gauss-Seidel (RGS) algorithm (Algorithm 4, also called the randomized coordinate descent algorithm).
\begin{table}[htp]
\begin{center}
\begin{tabular*}{115mm}{l}
\toprule {\bf Algorithm 4}. Randomized Gauss-Seidel \cite{leventhal2010rando}
\\ \hline \vspace{-3mm}\\ 
\qquad Initialize $\mbf x^0\in\mbbr^n$ \\
\qquad {\bf for} $k=1,2,\ldots$ {\bf do}\\
\qquad\qquad Pick $j\in[n]$ with probability $\|{\bf a}_j\|_2^2/\|{\bf A}\|_\rmf^2$\\
\qquad\qquad Set $\dsp{\bf x}^k={\bf x}^{k-1}-\frac{{\bf a}_j^\rmt({\bf Ax}^{k-1}-{\bf b})}{\|{\bf a}_j\|_2^2}{\bf e}_j$\\
\bottomrule
\end{tabular*}
\end{center}
\end{table}
The following theorem is a restatement of Lemma 4.2 of \cite{ma2015conve} and will be used to prove the refined bound for REGS. Here we provide a proof for completeness. 

\theorem\label{AAdag} Let ${\bf A}\in\mbbr^{m\times n}$ and ${\bf b}\in\mbbr^m$. Let ${\bf x}^k$ denote the $k$th iterate of {\rm RGS} applied to ${\bf Ax=b}$ with arbitrary ${\bf x}^0\in\mbbr^n$. In exact arithmetic, it holds $$\mbbe\bem\|{\bf Ax}^{k}-{\bf AA^\dag b}\|_2^2\eem\leq \rho^k\|{\bf Ax}^0-{\bf AA^\dag b}\|_2^2.$$
\endtheorem
\proof 
By ${\bf A^\rmt b}={\bf A^\rmt AA}^\dag{\bf b}$, we have \beqas {\bf x}^k-{\bf A^\dag b}&=&{\bf x}^{k-1}-{\bf A^\dag b}-\frac{{\bf a}_j^\rmt({\bf Ax}^{k-1}-{\bf b})}{\|{\bf a}_j\|^2_2}{\bf e}_j\\&=&{\bf x}^{k-1}-{\bf A^\dag b}- \frac{{\bf e}_j^\rmt({\bf A^\rmt Ax}^{k-1}-{\bf A^\rmt b})}{\|{\bf a}_j\|^2_2}{\bf e}_j\\&=&{\bf x}^{k-1}-{\bf A^\dag b}- \frac{{\bf e}_j^\rmt{\bf A^\rmt A}({\bf x}^{k-1}-{\bf A^\dag b})}{\|{\bf a}_j\|^2_2}{\bf e}_j\\&=&\l({\bf I}-\frac{{\bf e}_j{\bf e}_j^\rmt{\bf A^\rmt A}}{\|{\bf a}_j\|^2_2}\r)({\bf x}^{k-1}-{\bf A^\dag b}),\eeqas which yields $$ {\bf Ax}^k-{\bf AA^\dag b}=\l({\bf I}-\frac{{\bf a}_j{\bf a}_j^\rmt}{\|{\bf a}_j\|^2_2}\r)({\bf Ax}^{k-1}-{\bf AA^\dag b}).$$ 
It follows that \beqas&&\mbbe_{k-1}\bem\|{\bf Ax}^k-{\bf AA^\dag b}\|_2^2\eem=\mbbe_{k-1}\bem({\bf Ax}^k-{\bf AA^\dag b})^\rmt({\bf Ax}^k-{\bf AA^\dag b})\eem\\&=&\mbbe_{k-1}\bem({\bf Ax}^{k-1}-{\bf AA^\dag b})^\rmt\dsp\l({\bf I}-\frac{{\bf a}_j{\bf a}_j^\rmt}{\|{\bf a}_j\|^2_2}\r)^2({\bf Ax}^{k-1}-{\bf AA^\dag b})\eem\\&=&\mbbe_{k-1}\bem({\bf Ax}^{k-1}-{\bf AA^\dag b})^\rmt\dsp\l({\bf I}-\frac{{\bf a}_j{\bf a}_j^\rmt}{\|{\bf a}_j\|^2_2}\r)({\bf Ax}^{k-1}-{\bf AA^\dag b})\eem\\&=&({\bf Ax}^{k-1}-{\bf AA^\dag b})^\rmt\dsp\l({\bf I}-\frac{\bf AA^\rmt}{\|{\bf A}\|^2_\rmf}\r)({\bf Ax}^{k-1}-{\bf AA^\dag b})\\&\leq&\rho\|{\bf Ax}^{k-1}-{\bf AA^\dag b}\|_2^2.\quad (\mbox{by Lemma \ref{leq}}) \eeqas Taking expectation gives $$\mbbe\bem\|{\bf Ax}^k-{\bf AA^\dag b}\|_2^2\eem\leq\rho\mbbe\bem\|{\bf Ax}^{k-1}-{\bf AA^\dag b}\|_2^2\eem.$$ Unrolling the recurrence yields the result.
\endproof

If $\bf A$ has full column rank, Theorem \ref{AAdag} implies that $\mbf x^k$ converges linearly in expectation to $\bf A^\dag b$. If $\bf A$ does not have full column rank, RGS fails to converge (see \cite[section 3.3]{ma2015conve}). Ma, Needell, and Ramdas \cite{ma2015conve} proposed the following  randomized extended Gauss-Seidel algorithm (we call it REGS-MNR, see Algorithm 5) to resolve this problem. 
\begin{table}[htp]
\begin{center}
\begin{tabular*}{115mm}{l}
\toprule {\bf Algorithm 5}. REGS-MNR \cite{ma2015conve}
\\ \hline \vspace{-3mm}\\ 
\qquad Initialize $\mbf x^0\in\mbbr^n$ and $\mbf z^0\in\mbf x^0+\ran(\mbf A^\rmt)$ \\
\qquad {\bf for} $k=1,2,\ldots$ {\bf do}\\
\qquad\qquad Pick $j\in[n]$ with probability $\|{\bf a}_j\|_2^2/\|{\bf A}\|_\rmf^2$\\
\qquad\qquad Set $\dsp{\bf x}^k={\bf x}^{k-1}-\frac{{\bf a}_j^\rmt({\bf Ax}^{k-1}-{\bf b})}{\|{\bf a}_j\|_2^2}{\bf e}_j$\\
\qquad\qquad Pick $i\in[m]$ with probability $\|\wt{\bf a}_i\|_2^2/\|{\bf A}\|_\rmf^2$\\
\qquad\qquad Set $\dsp\mbf P_i=\mbf I-{\wt{\bf a}_i\wt{\bf a}_i^\rmt}/{\|\wt{\bf a}_i\|^2_2}$\\
\qquad\qquad Set $\dsp{\bf z}^k=\mbf P_i({\bf z}^{k-1}+\mbf x^k-\mbf x^{k-1})$\\
\qquad Output $\mbf x^t-\mbf z^t$ at some step $t$ as the estimated solution\\
\bottomrule
\end{tabular*}
\end{center}
\end{table}
\remark The original randomized extended Gauss-Seidel algorithm of {\rm\cite{zouzias2013rando}} uses $\mbf x^0=\mbf 0$ and ${\bf z}^0=\mbf 0$. Here, we use $\mbf x^0\in\mbbr^n$ and $\mbf z^0\in\mbf x^0+\ran(\mbf A^\rmt)$.
\endremark

Ma, Needell, and Ramdas proved that REGS-MNR converges linearly in expectation to the least norm solution in the case that $\mbf A$ has full row rank and $m<n$. They provided the convergence bound (see \cite[Theorem 4.1]{ma2015conve} for details) \beq\label{regsb}\mbbe\bem\|\mbf x^k-\mbf z^k-\mbf A^\dag\mbf b\|_2^2\eem\leq\rho^k\|\mbf A^\dag\mbf b\|_2^2+\frac{2\rho^{\lf k/2\rf}}{\sigma_r^2(\mbf A)}\|{\bf AA^\dag b}\|_2^2.\eeq 
Their proof (see \cite[Page 1600, lines 10-11]{ma2015conve}) uses Fact 1 of \cite[Page 1598]{ma2015conve} (which is that for any $\mbf P_i$ as in Algorithm 5, $\mbbe\bem\|\mbf P_i\mbf v\|_2^2\eem\leq\rho\|\mbf v\|_2^2$ for any vector $\mbf v\in\mbbr^n$) to show that  $$\mbbe\bem\|\mbf P_i(\mbf x^{k-1}-\mbf z^{k-1}-\mbf A^\dag\mbf b)\|_2^2\eem\leq\rho\|\mbf x^{k-1}-\mbf z^{k-1}-\mbf A^\dag\mbf b\|_2^2.$$ 
However, Fact 1 of \cite[Page 1598]{ma2015conve} is not true because for any nonzero  vector $\mbf v\in\nul(\bf A)$, we have $$\mbbe\bem\|\mbf P_i\mbf v\|_2^2\eem=\mbbe\bem\mbf v^\rmt\mbf P_i\mbf v\eem=\mbf v^\rmt\l(\mbf I-\frac{\bf A^\rmt A}{\|\mbf A\|_\rmf^2}\r)\mbf v=\|\mbf v\|_2^2.$$ 
 Therefore, the proof is incomplete. This issue can be resolved easily. Actually, by Lemma \ref{leq}, the above inequality still holds because $\mbf x^{k-1}-\mbf z^{k-1}-\mbf A^\dag\mbf b\in\ran(\mbf A^\rmt)$, which can be proved by induction.

\begin{table}[htp]
\begin{center}
\begin{tabular*}{115mm}{l}
\toprule {\bf Algorithm 6}. REGS-E  
\\ \hline \vspace{-3mm}\\ 
\qquad Initialize $\mbf x^0\in\mbbr^n$ and $\mbf z^0\in\ran(\mbf A^\rmt)$\\
\qquad {\bf for} $k=1,2,\ldots$ {\bf do}\\
\qquad\qquad Pick $j\in[n]$ with probability $\|{\bf a}_j\|_2^2/\|{\bf A}\|_\rmf^2$\\
\qquad\qquad Set $\dsp{\bf x}^k={\bf x}^{k-1}-\frac{{\bf a}_j^\rmt({\bf Ax}^{k-1}-{\bf b})}{\|{\bf a}_j\|_2^2}{\bf e}_j$\\
\qquad\qquad Pick $i\in[m]$ with probability $\|\wt{\bf a}_i\|_2^2/\|{\bf A}\|_\rmf^2$\\
\qquad\qquad Set $\dsp{\bf z}^k={\bf z}^{k-1}-\frac{\wt{\bf a}_i^\rmt({\bf z}^{k-1}-{\bf x}^k)}{\|\wt{\bf a}_i\|^2_2}\wt{\bf a}_i$\\
\bottomrule
\end{tabular*}
\end{center}
\end{table}

Next we study the convergence of REGS for a general linear system (consistent or inconsistent, full rank or rank-deficient). For the convenience of discussion, we present the following randomized extended Gauss-Seidel algorithm (we call it REGS-E, see Algorithm 6) which is mathematically equivalent to REGS-MNR. Actually, in exact  arithmetic, the vector $\mbf z^k$ in REGS-E is equal to the vector $\mbf x^k-\mbf z^k$ in REGS-MNR.

\remark\label{RGSRK}  Essentially, {\rm REGS-E} is an {\rm RGS-RK} approach: $\mbf x^k$ is the $k$th iterate of {\rm RGS} and $\mbf z^k$ is one step {\rm RK} update for the linear system ${\bf Az=Ax}^k$ from $\mbf z^{k-1}$.
\endremark

In the following theorem, we show that the vector ${\bf z}^k$ in REGS-E converges linearly in expectation to $\bf A^\dag b$. Our proof is almost the same as that of \cite[Theorem 4.1]{ma2015conve} but avoiding their mistake. The new convergence bound is attainable (see Remark \ref{sharp2}) and  obviously better than the bound $(\ref{regsb})$. 

\theorem\label{main2} Let ${\bf A}\in\mbbr^{m\times n}$ and ${\bf b}\in\mbbr^m$. Let ${\bf z}^k$ denote the $k$th iterate of {\rm REGS-E} with arbitrary $\mbf x^0\in\mbbr^n$ and $\mbf z^0\in\ran(\mbf A^\rmt)$. In exact arithmetic, it holds  \beq\label{bound3}\mbbe\bem\|{\bf z}^k-{\bf A^\dag b}\|_2^2\eem\leq\rho^k\|{\bf z}^0-{\bf A^\dag b}\|_2^2+\frac{\rho^k(1-\rho^k)}{\sigma_{r}^2(\mbf A)}\|{\bf Ax}^0-{\bf AA^\dag b}\|_2^2.\eeq
\endtheorem
\proof By $\mbf z^0\in\ran(\mbf A^\rmt)$ and $\mbf A^\dag\mbf b\in\ran(\mbf A^\rmt)$, we have $\mbf z^0-\mbf A^\dag\mbf b\in\ran(\mbf A^\rmt)$. Then it is easy to show that $\mbf z^k-\mbf A^\dag\mbf b\in\ran(\mbf A^\rmt)$ by induction.
We now analyze the norm of $\mbf z^k-{\bf A^\dag b}$. Note that \beqas\mbf z^k-{\bf A^\dag b}&=&{\bf z}^{k-1}-\frac{\wt{\bf a}_i^\rmt({\bf z}^{k-1}-{\bf x}^k)}{\|\wt{\bf a}_i\|^2_2}\wt{\bf a}_i-{\bf A^\dag b}\\&=&\dsp\l(\mbf I-\frac{\wt{\bf a}_i\wt{\bf a}_i^\rmt}{\|\wt{\bf a}_i\|_2^2}\r){\bf z}^{k-1}+\frac{\wt{\bf a}_i\wt{\bf a}_i^\rmt}{\|\wt{\bf a}_i\|_2^2}\mbf x^k-{\bf A^\dag b}\\&=&\dsp\l(\mbf I-\frac{\wt{\bf a}_i\wt{\bf a}_i^\rmt}{\|\wt{\bf a}_i\|_2^2}\r)({\bf z}^{k-1}-{\bf A^\dag b})+\frac{\wt{\bf a}_i\wt{\bf a}_i^\rmt}{\|\wt{\bf a}_i\|_2^2}(\mbf x^k-{\bf A^\dag b}).\eeqas It follows from the orthogonality, i.e., $$(\mbf x^k-{\bf A^\dag b})^\rmt\frac{\wt{\bf a}_i\wt{\bf a}_i^\rmt}{\|\wt{\bf a}_i\|_2^2}\l(\mbf I-\frac{\wt{\bf a}_i\wt{\bf a}_i^\rmt}{\|\wt{\bf a}_i\|_2^2}\r)({\bf z}^{k-1}-{\bf A^\dag b})=0,$$ that \beqa\label{sum}\|\mbf z^k-{\bf A^\dag b}\|_2^2&=&\l\|\dsp\l(\mbf I-\frac{\wt{\bf a}_i\wt{\bf a}_i^\rmt}{\|\wt{\bf a}_i\|_2^2}\r)({\bf z}^{k-1}-{\bf A^\dag b})\r\|_2^2+\l\|\frac{\wt{\bf a}_i\wt{\bf a}_i^\rmt}{\|\wt{\bf a}_i\|_2^2}(\mbf x^k-{\bf A^\dag b})\r\|_2^2.
\eeqa 
It follows from \beqas&&\mbbe_{k-1}\bem\l\|\dsp\l(\mbf I-\frac{\wt{\bf a}_i\wt{\bf a}_i^\rmt}{\|\wt{\bf a}_i\|_2^2}\r)({\bf z}^{k-1}-{\bf A^\dag b})\r\|_2^2\eem\\&=&\mbbe_{k-1}\bem({\bf z}^{k-1}-{\bf A^\dag b})^\rmt\dsp\l(\mbf I-\frac{\wt{\bf a}_i\wt{\bf a}_i^\rmt}{\|\wt{\bf a}_i\|_2^2}\r)^2({\bf z}^{k-1}-{\bf A^\dag b})\eem\\&=&\mbbe_{k-1}\bem({\bf z}^{k-1}-{\bf A^\dag b})^\rmt\dsp\l(\mbf I-\frac{\wt{\bf a}_i\wt{\bf a}_i^\rmt}{\|\wt{\bf a}_i\|_2^2}\r)({\bf z}^{k-1}-{\bf A^\dag b})\eem
\\&=&({\bf z}^{k-1}-{\bf A^\dag b})^\rmt\dsp\l(\mbf I-\frac{\bf A^\rmt A}{\|\mbf A\|_\rmf^2}\r)({\bf z}^{k-1}-{\bf A^\dag b})\\&\leq&\rho\|{\bf z}^{k-1}-{\bf A^\dag b}\|_2^2\quad (\mbox{by Lemma \ref{leq}})\eeqas that \beq\label{sum1}\mbbe\bem\l\|\dsp\l(\mbf I-\frac{\wt{\bf a}_i\wt{\bf a}_i^\rmt}{\|\wt{\bf a}_i\|_2^2}\r)({\bf z}^{k-1}-{\bf A^\dag b})\r\|_2^2\eem\leq\rho\mbbe\bem\|{\bf z}^{k-1}-{\bf A^\dag b}\|_2^2\eem.\eeq 
It follows from \beqas&&\mbbe_{k-1}\bem\dsp\l\|\frac{\wt{\bf a}_i\wt{\bf a}_i^\rmt}{\|\wt{\bf a}_i\|_2^2}(\mbf x^k-{\bf A^\dag b})\r\|_2^2\eem\\&=&\mbbe_{k-1}\bem\dsp(\mbf x^k-{\bf A^\dag b})^\rmt\l(\frac{\wt{\bf a}_i\wt{\bf a}_i^\rmt}{\|\wt{\bf a}_i\|_2^2}\r)^2(\mbf x^k-{\bf A^\dag b})\eem\\&=&\mbbe_{k-1}^j\bem\mbbe_{k-1}^i\bem\dsp(\mbf x^k-{\bf A^\dag b})^\rmt\frac{\wt{\bf a}_i\wt{\bf a}_i^\rmt}{\|\wt{\bf a}_i\|_2^2}(\mbf x^k-{\bf A^\dag b})\eem\eem
\\&=&\mbbe_{k-1}^j\bem\dsp(\mbf x^k-{\bf A^\dag b})^\rmt\frac{\bf A^\rmt A}{\|\mbf A\|_\rmf^2}(\mbf x^k-{\bf A^\dag b})\eem\\&=&\frac{1}{\|\mbf A\|_\rmf^2}\mbbe_{k-1}\bem\|{\bf Ax}^k-{\bf AA^\dag b}\|_2^2\eem\eeqas that \beqa\nn\mbbe\bem\dsp\l\|\frac{\wt{\bf a}_i\wt{\bf a}_i^\rmt}{\|\wt{\bf a}_i\|_2^2}(\mbf x^k-{\bf A^\dag b})\r\|_2^2\eem&=&\frac{1}{\|\mbf A\|_\rmf^2}\mbbe\bem\|{\bf Ax}^k-{\bf AA^\dag b}\|_2^2\eem\\\label{sum2}&\leq&\frac{\rho^k}{\|\mbf A\|_\rmf^2}\|{\bf Ax}^0-{\bf AA^\dag b}\|_2^2. \quad(\mbox{by Theorem \ref{AAdag}})\eeqa
Combining (\ref{sum}), (\ref{sum1}), and (\ref{sum2}) yields \beqas\mbbe\bem\|\mbf z^k-{\bf A^\dag b}\|_2^2\eem&=&\mbbe\bem\l\|\dsp\l(\mbf I-\frac{\wt{\bf a}_i\wt{\bf a}_i^\rmt}{\|\wt{\bf a}_i\|_2^2}\r)({\bf z}^{k-1}-{\bf A^\dag b})\r\|_2^2\eem+\mbbe\bem\l\|\dsp\frac{\wt{\bf a}_i\wt{\bf a}_i^\rmt}{\|\wt{\bf a}_i\|_2^2}(\mbf x^k-{\bf A^\dag b})\r\|_2^2\eem\\&\leq&\rho\mbbe\bem\|{\bf z}^{k-1}-{\bf A^\dag b}\|_2^2\eem+\frac{\rho^k}{\|\mbf A\|_\rmf^2}\|{\bf Ax}^0-{\bf AA^\dag b}\|_2^2\\&\leq&\cdots\leq\rho^k\|{\bf z}^0-{\bf A^\dag b}\|_2^2+\frac{\rho^k}{\|\mbf A\|_\rmf^2}\|{\bf Ax}^0-{\bf AA^\dag b}\|_2^2\sum_{l=0}^{k-1}\rho^l\\&=&\rho^k\|{\bf z}^0-{\bf A^\dag b}\|_2^2+\frac{\rho^k}{\|\mbf A\|_\rmf^2}\|{\bf Ax}^0-{\bf AA^\dag b}\|_2^2\frac{1-\rho^k}{1-\rho}\\&=&\rho^k\|{\bf z}^0-{\bf A^\dag b}\|_2^2+\frac{\rho^k(1-\rho^k)}{\sigma_{r}^2(\mbf A)}\|{\bf Ax}^0-{\bf AA^\dag b}\|_2^2.\eeqas This completes the proof.
\endproof

\remark\label{sharp2} By Lemma {\rm 1}, if $\sigma_1(\mbf A)=\sigma_r(\mbf A)$, then all the inequalities in the proofs of Theorems {\rm \ref{AAdag}} and {\rm \ref{main2}} become equalities. 
\endremark

\section{Numerical results} We compare the bounds (\ref{rekb}), (\ref{bound2}),  (\ref{regsb}), and (\ref{bound3}) via a set of small examples. For given $m$, $n$,  $r=\rank(\mbf A)$, $\sigma_1(\mbf A)$, and $\sigma_r(\mbf A)$, we construct a matrix $\bf A$ by $\bf A=UD V^\rmt$, where $\mbf U\in\mbbr^{m\times r}$ and $\mbf V\in\mbbr^{n\times r}$. Entries of $\mbf U$ and $\mbf V$ are generated from a standard normal distribution and then columns are orthonormalized. The matrix $\mbf D$ is an $r\times r$ diagonal matrix whose first $r-2$ diagonal entries are uniformly distributed numbers in $[\sigma_r(\mbf A), \sigma_1(\mbf A)]$, and the last two diagonal entries are $\sigma_r(\mbf A)$ and $\sigma_1(\mbf A)$. 

We consider four cases: (i) $\bf Ax=b$ is consistent and $\rank(\mbf A)=n$; (ii) $\bf Ax=b$ is consistent and $\rank(\mbf A)<n$; (iii) $\bf Ax=b$ is inconsistent and $\rank(\mbf A)=n$; (iv) $\bf Ax=b$ is inconsistent and $\rank(\mbf A)<n$. To construct a consistent linear system, we set $\bf b=Ax$ where $\bf x$ is a vector with entries generated from a standard normal distribution. 
To construct an inconsistent linear system, we set $\bf b=Ax+r$ where $\bf x$ is a vector with entries generated from a standard normal distribution and the residual $\mbf r\in\nul(\mbf A^\rmt)$. Note that one can obtain such a vector $\mbf r$ by the MATLAB function {\tt null}. 
\begin{figure}[!htpb]
\centerline{\epsfig{figure=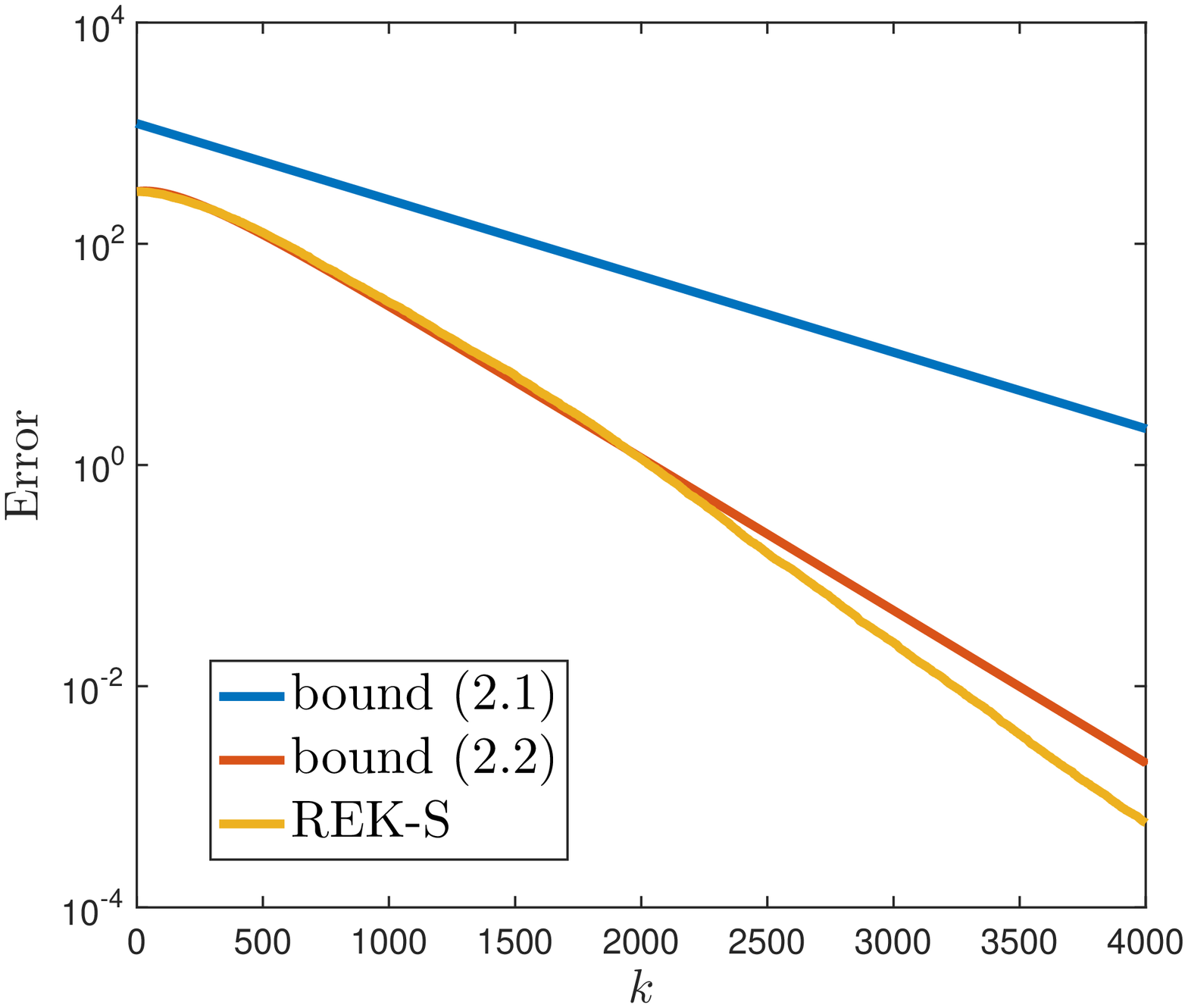,height=2.1in}\epsfig{figure=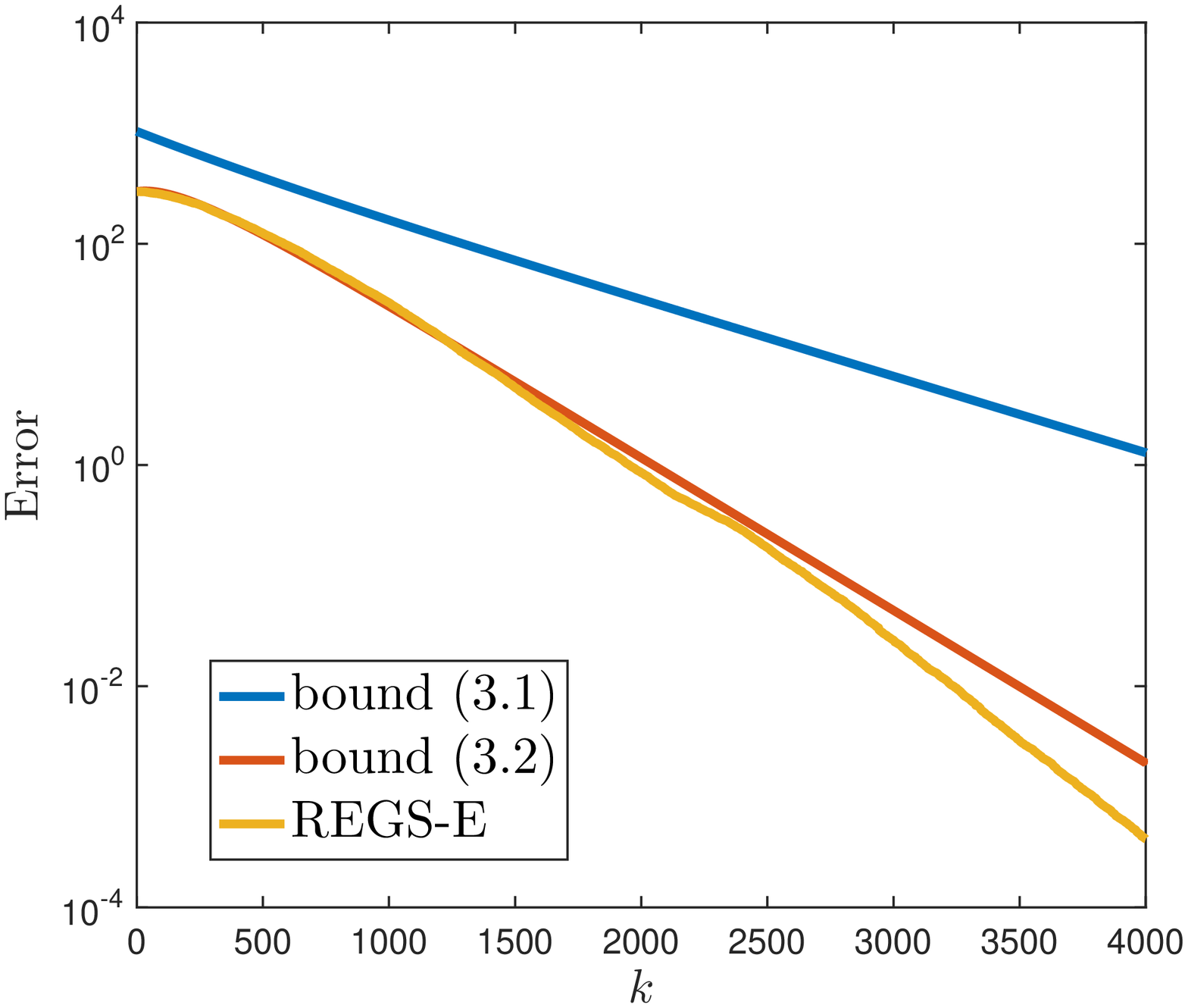,height=2.1in}}
\caption{The error $\|\mbf x^k-\mbf A^\dag\mbf b\|_2^2$ for {\rm REK-S (}left{\rm)} and the error $\|\mbf z^k-\mbf A^\dag\mbf b\|_2^2$ for {\rm REGS-E (}right{\rm)} on a consistent linear system with full column rank $\bf A$: $m=500$, $n=250$, $r=250$, $\sigma_1(\mbf A)=1.25$, and $\sigma_r(\mbf A)=1$.} \label{fig1}
\end{figure} 

\begin{figure}[!htpb]
\centerline{\epsfig{figure=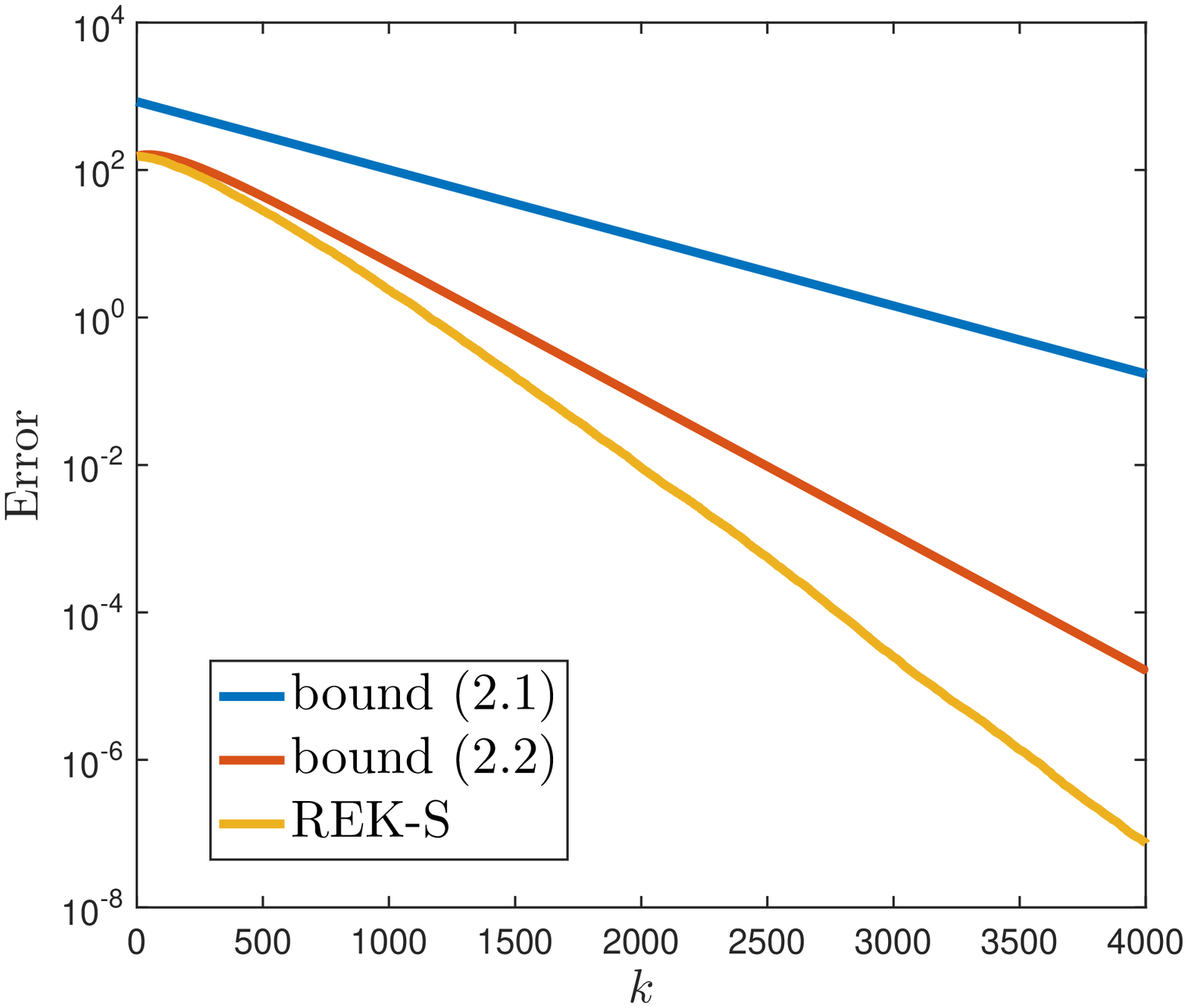,height=2.1in}\epsfig{figure=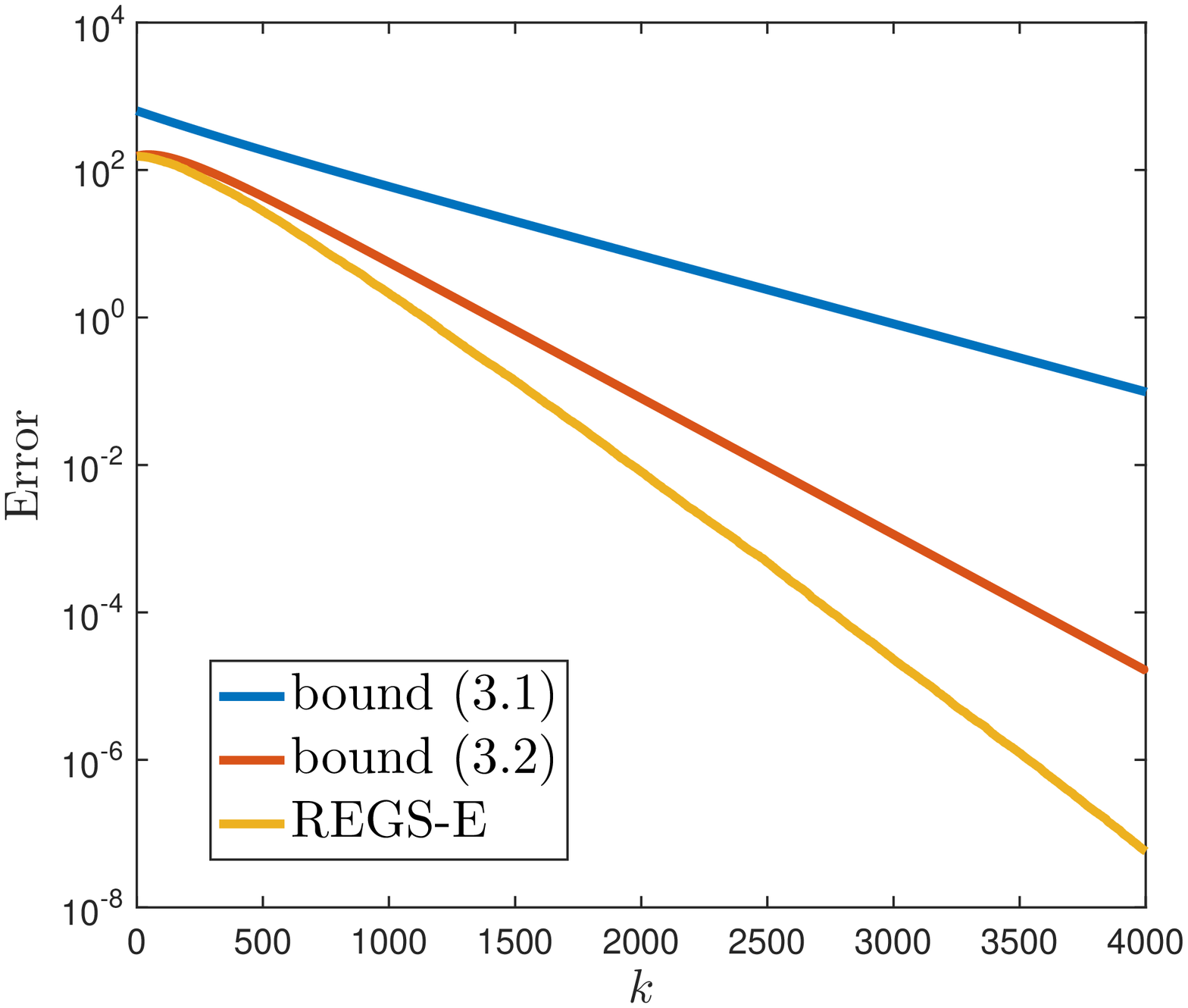,height=2.1in}}
\caption{The error $\|\mbf x^k-\mbf A^\dag\mbf b\|_2^2$ for {\rm REK-S (}left{\rm)} and the error $\|\mbf z^k-\mbf A^\dag\mbf b\|_2^2$ for {\rm REGS-E (}right{\rm)} on a consistent linear system with rank-deficient $\bf A$: $m=500$, $n=250$, $r=150$, $\sigma_1(\mbf A)=1.5$, and $\sigma_r(\mbf A)=1$.}\label{fig2}
\end{figure} 

\begin{figure}[!htpb]
\centerline{\epsfig{figure=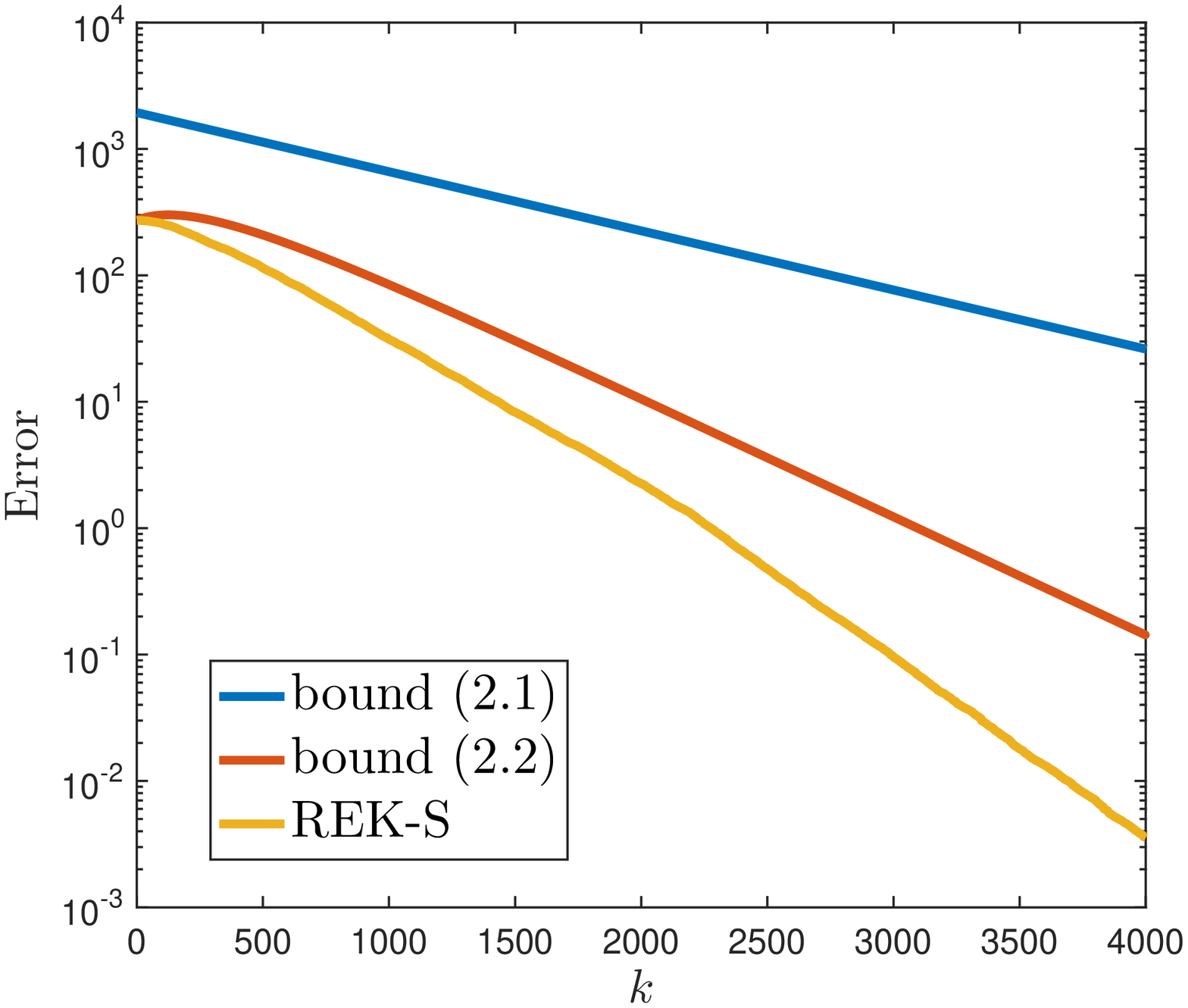,height=2.1in}\epsfig{figure=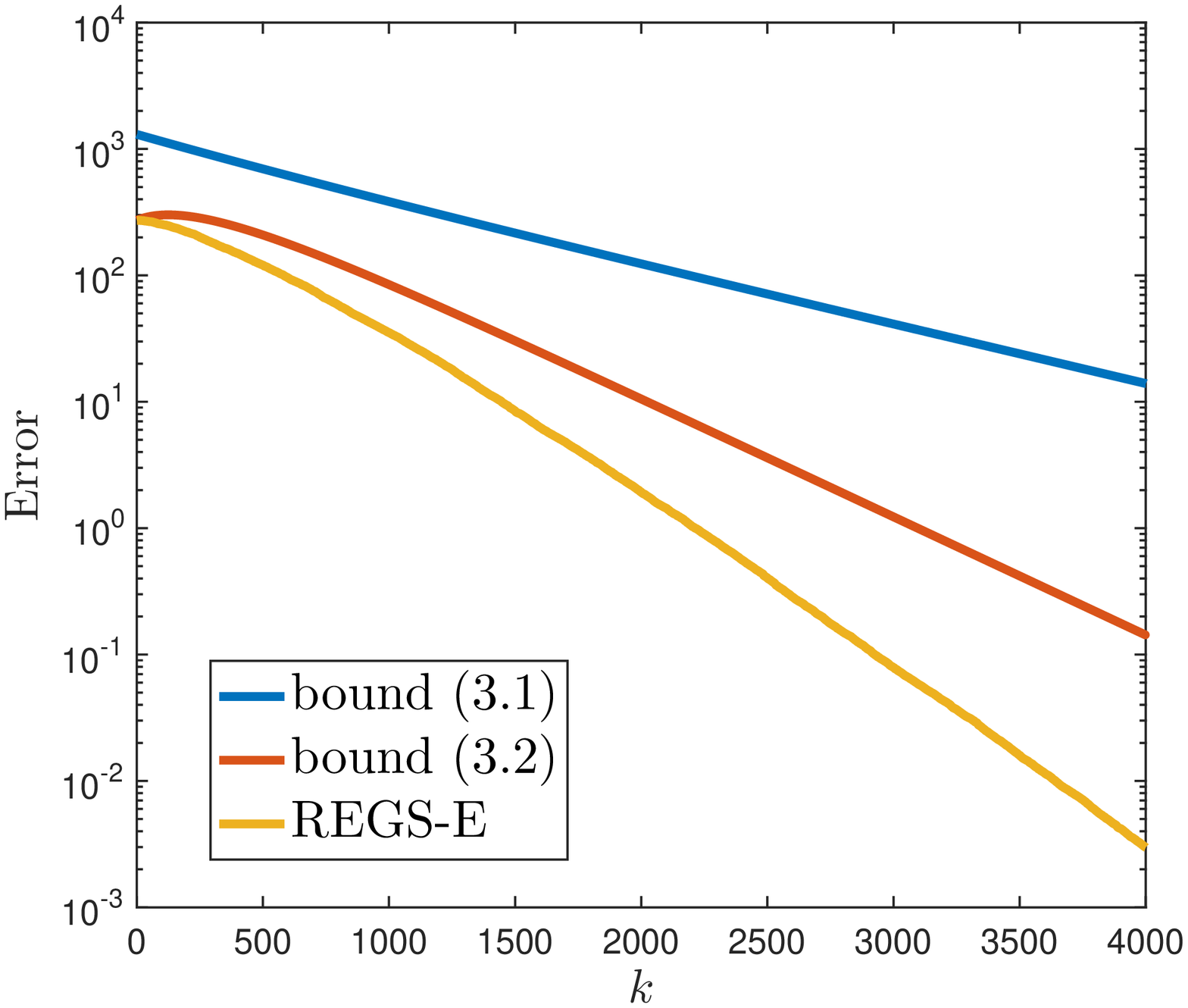,height=2.1in}}
\caption{The error $\|\mbf x^k-\mbf A^\dag\mbf b\|_2^2$ for {\rm REK-S (}left{\rm)} and the error $\|\mbf z^k-\mbf A^\dag\mbf b\|_2^2$ for {\rm REGS-E (}right{\rm)} on an inconsistent linear system with full column rank $\bf A$: $m=500$, $n=250$, $r=250$, $\sigma_1(\mbf A)=1.75$, and $\sigma_r(\mbf A)=1$.}\label{fig3}
\end{figure} 

\begin{figure}[!htpb]
\centerline{\epsfig{figure=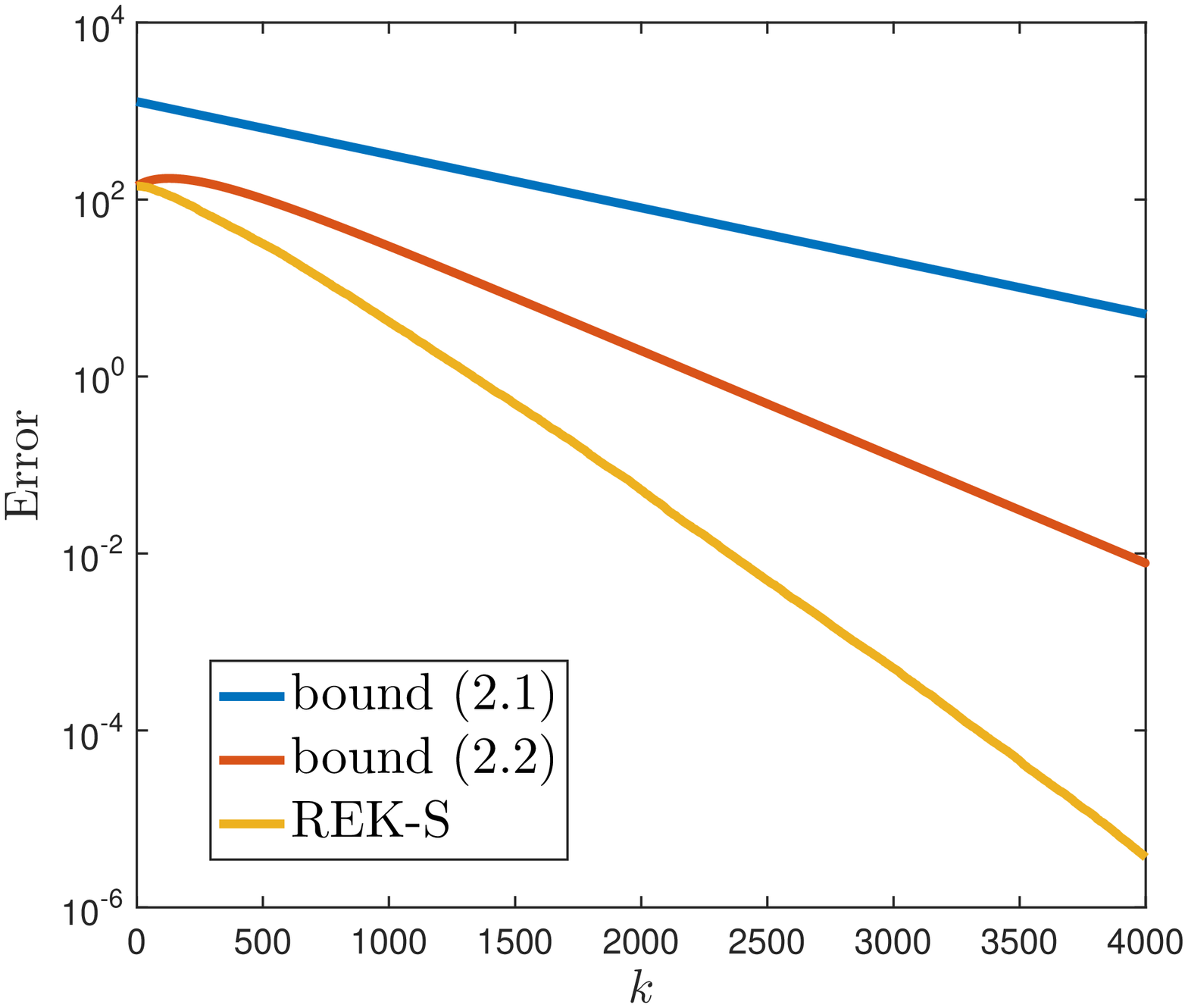,height=2.1in}\epsfig{figure=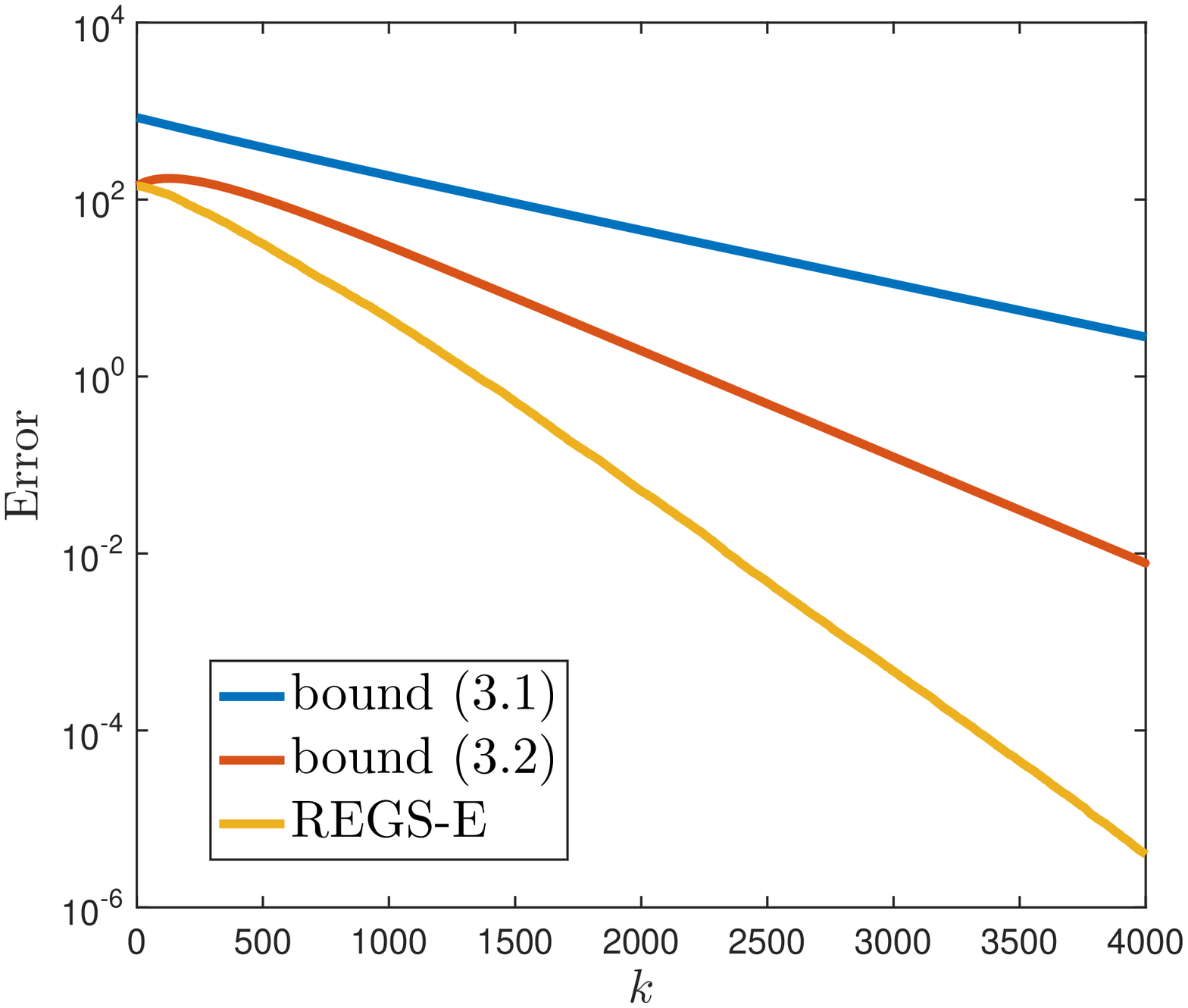,height=2.1in}}
\caption{The error $\|\mbf x^k-\mbf A^\dag\mbf b\|_2^2$ for {\rm REK-S (}left{\rm)} and the error $\|\mbf z^k-\mbf A^\dag\mbf b\|_2^2$ for {\rm REGS-E (}right{\rm)} on an inconsistent linear system with rank-deficient $\bf A$: $m=500$, $n=250$, $r=150$, $\sigma_1(\mbf A)=2$, and $\sigma_r(\mbf A)=1$.}\label{fig4}
\end{figure} 

In Figures \ref{fig1}-\ref{fig4}, we plot the error $\|\mbf x^k-\mbf A^\dag\mbf b\|_2^2$ for REK-S with $\mbf x^0=\mbf 0$ and $\mbf z^0=\mbf b$ and the error $\|\mbf z^k-\mbf A^\dag\mbf b\|_2^2$ for REGS-E with $\mbf x^0=\mbf 0$ and $\mbf z^0=\mbf 0$ for the four cases, respectively. For each case, we average the error over 20 trials for the same problem. For all cases, our bounds (\ref{bound2}) and (\ref{bound3}) are much better than the existing bounds (\ref{rekb}) and (\ref{regsb}).

\section{Conclusion} We have shown that REK is essentially an RK-RK approach and that REGS is essentially an RGS-RK approach. We have proposed refined upper bounds for the convergence of both algorithms. These upper bounds are attained for the case that all nonzero  singular values of $\bf A$ are the same. Our convergence analysis applies to all types of linear systems. 
The acceleration technique such as that used in \cite{liu2016accel} is being considered.


\end{document}